\documentclass[11pt]{article}
\usepackage{amsmath}
\usepackage{amsfonts}
\usepackage[bookmarksnumbered=true]{hyperref}

\newtheorem{theorem}{Theorem}[section]

\newtheorem{lemma}[theorem]{Lemma}

\newtheorem{remark}[theorem]{Remark}

\begin{document}

\title{The Equivalence between Uniqueness and Continuous Dependence
 of Solution for BDSDEs\thanks{
 This work is supported by National Natural Science Foundation of China Grant
10771122, Natural Science Foundation of Shandong Province of China
Grant Y2006A08 and National Basic Research Program of China (973
Program, No.2007CB814900)}}
\author{Qingfeng Zhu$^{\rm a}$ and Yufeng Shi$^{\rm b}$\thanks{Corresponding author, E-mail: yfshi@sdu.edu.cn}\\
{\small $^{\rm a}$ School of Statistics and Mathematics, Shandong
University of Finance},\\
{\small Jinan 250014, China}\\
{\small$^{\rm b}$School of Mathematics, Shandong University, Jinan 250100, China}}
\maketitle

\begin{abstract}In this paper, we
prove that, if the coefficient $f=f(t,y,z)$ of backward doubly
stochastic differential equations (BDSDEs for short) is assumed to
be continuous and linear growth in $(y,z),$ then the uniqueness of
solution and continuous dependence with respect to the coefficients
$f$, $g$ and the terminal value $\xi$ are equivalent.\\
\indent{\it keywords:} backward doubly stochastic differential equations, uniqueness,
continuous dependence
\end{abstract}

\section{Introduction}\label{sec:1}

Nonlinear backward stochastic differential equations (BSDEs in short)
have been independently introduced by Pardoux and Peng \cite{PP1} and Duffie
and Epstein \cite{DE}. Since then, BSDEs have been studied intensively. In
particular, many efforts have been made to relax the assumption on the
generator. For instance, Lepeltier and San Martin \cite{LS} have proved the
existence of a solution for the case when the generator is only continuous
with linear growth, and Jia and Peng \cite{JP} obtained  that BSDE  has either
one or uncountably many solutions, if the generator satisfies the conditions
given in \cite{LS}. Jia and Yu \cite{JY} studied the equivalence between uniqueness
and continuous dependence of solution for BSDEs with continuous coefficient.
Another main reason is due to their enormous range of
applications in such diverse fields as mathematical finance (see \cite{DE} and El
Karoui et al. \cite{EPQ}, partial differential equations (see Peng \cite{P}),
stochastic control (see Ji and Wu \cite{JW}), nonlinear mathematical expectations
(see Jiang \cite{J} and Fan \cite{F}), and so on.

A class of backward doubly stochastic differential equations (BDSDEs
in short) was introduced by Pardoux and Peng \cite{PP2} in 1994, in order
to provide a probabilistic interpretation for the solutions of a
class of semilinear stochastic partial differential equations (SPDEs
in short). They have proved the existence and uniqueness of
solutions for BDSDEs under uniformly Lipschitz conditions. Since
then, Shi et al. \cite{SGL} have relaxed the Lipschitz assumptions to
linear growth conditions. Bally and Matoussi \cite{BM} have given a
probabilistic interpretation of the solutions in Sobolev spaces for
semilinear parabolic SPDEs in terms of BDSDEs. Zhang and Zhao \cite{ZZ}
have proved the existence and uniqueness of solution for BDSDEs on
infinite horizons, and described the stationary solutions of SPDEs
by virtue of the solutions of BDSDEs on infinite horizons. Recently,
Ren et al. \cite{RLH} and Hu and Ren \cite{HR} considered the
BDSDEs driven by Levy process with Lipschitz coefficient and applications in SPDEs

Because of their important significance to SPDEs, it is necessary to
give intensive investigation to the theory of BDSDEs.
In this paper we will prove that if the coefficient $f$ satisfis the
conditions given in \cite{SGL}, then the uniqueness of solution and
continuous dependence with respect to $f$, $g$ and $\xi$ are
equivalent. We consider the following 1-dimesional  backward
doubly stochastic differential equations:
\begin{equation}\label{eq:1}
\begin{array}{lll}
Y_t=\xi+\int_t^Tf(s,Y_s,Z_s)ds+\int_t^Tg(s,Y_s,Z_s)dB_s-
\int_t^TZ_sdW_s,\ 0\leq t\leq T,
\end{array}
\end{equation}
where $\{W_t;0\leq t\leq T\}$ and $\{B_t;0\leq t\leq T\}$ are two
mutually independent standard Brownian Motions with values in $\mathbb{R}^d$
and $\mathbb{R}^l$, respectively, defined on $(\Omega,{\cal F},P)$. The
terminal condition $\xi$ and the coefficients $f = f(t, y, z)$ and
$g = g(t, y, z)$ are given.
 The solution $(Y_t, Z_t)_{t\in [0,T]}$ is a pair of square
integrable processes. An interesting problem is: what is the
relationship between uniqueness of solution and continuous
dependence with respect to $f$, $g$ and $\xi$? In the standard
situation where $f$ satisfies Lipschitz condition in $(y, z)$, it
was proved by Pardoux and Peng \cite{PP2} that there exists a unique
solution. In this case, the continuous dependence with respect to
$f$ and $\xi$ is an obvious result. However in the case where $f$ is
only continuous in ($y, z)$, in place of the Lipschitz condition,
Shi et al. \cite{SGL} have proved that there is at least one solution.
In fact, there is either one or uncountable many solutions in this
situation (see Shi and Zhu \cite{SZ}). Does the uniqueness of solution of
BDSDEs also imply the continuous dependence with respect to $f$, $g$
and $\xi$?

This paper is organized as follows. In Section \ref{sec:2} we formulate the
problem accurately and give some preliminary results. Section \ref{sec:3} is
devoted to proving the equivalence of uniqueness and continuous
dependence with respect to terminal value $\xi$. Finally, in Section
\ref{sec:4} we will prove the equivalence of uniqueness and continuous
dependence with respect to parameters $f$, $g$ and $\xi$.

\section{Preliminary}\label{sec:2}

{\bf Notation} The Euclidean norm of a vector
$x\in{\mathbb{R}}^k$ will be denoted by $|x|$, and for a
$d\times k$ matrix $A$, we define $|A|=\sqrt{TrAA^*}$, where $A^*$ is the transpose of $A$.

Let $(\Omega,{\cal{F}},P)$ be a probability space, and $T>0$ be an
arbitrarily fixed constant throughout this paper. Let $\{W_t;0\leq
t\leq T\}$ and $\{B_t;0\leq t\leq T\}$ be two mutually independent
standard Brownian Motions with values in $\mathbb{R}^d$ and $\mathbb{R}^l$,
respectively, defined on $(\Omega,\cal{F},P)$. Let $\cal{N}$ denote
the class of $P$-null sets of $\cal{F}$. For each $t\in [0,T]$, we
define ${\cal{F}}_t : ={\cal{F}}_t^W\vee {\cal{F}}_{t,T}^B$, where
for any process $\{\eta_t\}$,
${\cal{F}}_{s,t}^{\eta}=\sigma\{\eta_r-\eta_s;s\leq r\leq t\}\vee
{\cal{N}}$, ${\cal{F}}_t^{\eta}={\cal{F}}_{0,t}^{\eta}$.
 Note that the collection $\{{\cal{F}}_t;t\in [0,T]\}$ is neither
increasing nor decreasing, so it does not constitute a filtration.

We introduce the following notations:
\begin{eqnarray*}
S^2\left([0,T];\mathbb{R}^n\right)&=&\{v_t,0\leq t\leq T,\
 \mbox {is an}\ \mathbb{R}^n\mbox{-valued},\ {\cal F}_t \mbox{-measurable process}\\
&&
\mbox{such that} \ E(\sup_{0\leq t\leq T}|v_{t}|^{2})<\infty\},\\
M^2(0,T;\mathbb{R}^n)&=&\{v_t,0\leq t\leq T,\
 \mbox {is an}\ \mathbb{R}^n\mbox{-valued},\ {\cal F}_t\mbox{-measurable process}\\
&&
\mbox{such that} \ E\int_{0}^T|v_{t}|^{2}dt<\infty\}.
\end{eqnarray*}

Let
\begin{eqnarray*}
f:\Omega\times [0,T]\times \mathbb{R}\times \mathbb{R}^{ d}\rightarrow \mathbb{R},\quad
g:\Omega\times [0,T]\times \mathbb{R}\times \mathbb{R}^{d}\rightarrow \mathbb{R}^{l},
\end{eqnarray*}
 be jointly measurable such that
for any $(y,z)\in \mathbb{R}\times \mathbb{R}^{d}$,
\begin{eqnarray*}
f(\cdot,y,z)\in M^2(0,T;\mathbb{R}),\quad g(\cdot,y,z)\in
M^2(0,T;\mathbb{R}^{l}).
\end{eqnarray*}
and satisfy the following conditions:

\vspace{0.1cm}(H1) linear growth: $\exists\ 0<K< \infty$, such that

\vspace{0.2cm}\hspace{1cm}
$|f(\omega,t,y,z)|\leq K(1+|y|+|z|), \quad\forall\ ( \omega, t, y, z)\in
\Omega \times [0,T]\times \mathbb{R}\times \mathbb{R}^d;$

\vspace{0.2cm}(H2) For fixed $\omega$ and $t$, $f(\omega, t, \cdot,\cdot)$ is
continuous;

(H3) there exist constants $c>0$ and $0<\alpha<1$ such that

\vspace{0.2cm}\hspace{2cm}$|g(\omega,t,y^1,z^1)-g(\omega,t,y^2,z^2)|^2
\leq c|y^1-y^2|^2+\alpha |z^1-z^2|^2,$

\vspace{0.2cm}\noindent for all $(\omega,t)\in \Omega\times [0,T],
\ (y^1,z^1)\in \mathbb{R}\times \mathbb{R}^d, \ (y^2,z^2)\in \mathbb{R}\times \mathbb{R}^d.$

\begin{remark}\label{rmk:2.1}
In fact, (H1) can be replaced by the following condition:

(H4) there exist a constant $0<K< \infty$, such that

\vspace{0.2cm}
$|f(\omega,t,y,z)-f(\omega,0,0,0)|\leq K(1+|y|+|z|), \quad\forall\ ( \omega, t, y, z)\in
\Omega \times [0,T]\times \mathbb{R}\times \mathbb{R}^d.$

\end{remark}

In the sequel, it is not hard to check that all results in this paper also hold under
Assumptions (H2)-(H4).

By Theorem 4.1 in \cite{SGL}, under (H1)-(H3) and for each given
$\xi \in L^2(\Omega,{\cal{F}}_T,P)$, there exists at least one
solution $(Y_t, Z_t)_{t\in[0,T ]}\in S^2 \times M^2$ of BDSDE (\ref{eq:1}).
\cite{SGL} gives also the existence of the minimal solution
$(\underline{Y}_t,\underline{Z}_t)_{t\in [0,T]}$ of BDSDE (\ref{eq:1}) and
\cite{SZ} gives the maximal solution
$(\overline{Y}_t,\overline{Z}_t)_{t\in [0,T]}$ of BDSDE (\ref{eq:1}) in the
sense that any solution $(Y_t, Z_t)_{t\in[0,T]}\in S^2\times M^2$ of
BDSDE (\ref{eq:1}) must satisfy $\underline{Y}_t\leq Y_t\leq \overline{Y}_t,$
a.s., for all $t \in [0, T]$.

It is well known that under the standard assumptions where $f$ is
Lipschitz continuous in $(y, z)$, for any random variable $\xi$ in
$L^2(\Omega,{\cal{F}}_T,P)$, the BDSDE (\ref{eq:1}) has a unique adapted
solution, say $(Y_t, Z_t)_{t\in[0,T]}$ such that $Y \in S^2$ and
$Z\in M^2$ (see \cite{PP2}). And we have the following estimate for
solution of BDSDEs with Lipschitz continuous generator $f$ comes
from \cite{PP2}.

\begin{lemma}\label{lem:2.2} If $\xi^1,
\xi^2 \in L^2(\Omega,{\cal{F}}_T,P),$ $f$ is Lipschitz continuous in
$(Y, Z)$ and $g$ satisfies (H3). Then, for the solutions $(Y^1_t,
Z^1_t)_{t\in[0,T]}$ and $(Y^2_t, Z^2_t)_{t\in[0,T]}$ of the BDSDEs
$(f, g, T, \xi^1)$ and $(f,$ $g, T, \xi^2)$ respectively, we have
$$E[\sup\limits_{0\leq t\leq T}|Y^1_t-Y^2_t|^2]\leq C E|\xi^1-\xi^2|^2,$$
where $C$ is a positive constant only depending on Lipschitz
constants of $f$ and $g$.
\end{lemma}

Now, we recall some properties and associated approximation about BDSDEs
with $f$ and $g$ satisfying Assumptions (H1)-(H3) (see \cite{PP2} for details).

\begin{lemma}\label{lem:2.3} If $f$
satisfies Assumptions (H1) and (H2), and we set
\begin{eqnarray*}
\underline{f}_m(\omega,t,y,z)=\inf_{(y', z')\in Q^{1+d}}\{f(\omega,t,y',z')+m(|y-y'|+|z-z'|)\} 
\end{eqnarray*}
and
\begin{eqnarray*}
\overline{f}_m(\omega,t,y,z)=\sup_{(y', z')\in Q^{1+d}}\{f(\omega,t,y',z')-m(|y-y'|+|z-z'|)\} 
\end{eqnarray*}
then for any $m\geq K$, we have

{\rm (i)} linear growth: $\forall\ ( y, z)\in \mathbb{R}\times \mathbb{R}^d$ and $ t\in
[0,T],$ $$|\underline{f}_m(t,y,z)|\leq K(1+|y|+|z|),\ \mbox{and}\
|\overline{f}_m(t,y,z)|\leq K(1+|y|+|z|).$$

{\rm (ii)} monotonicity in $m$: $\forall\ ( y, z)\in \mathbb{R}\times \mathbb{R}^d$ and
$t\in [0,T],$ $\underline{f}_m(t,y,z)$ is non-decreasing in $m$ and
$\overline{f}_m(t,y,z)$ is non-increasing in $m$.

{\rm (iii)} Lipschitz condition: $\forall\ y_1, y_2 \in \mathbb{R}, z_1,z_1\in
\mathbb{R}^d$ and $t\in [0,T],$
$$|\underline{f}_m(t,y,z)-\underline{f}_m(t,y',z')|\leq
m(|y-y'|+|z-z'|),$$ and
$$|\overline{f}_m(t,y,z)-\overline{f}_m(t,y',z')|\leq
m(|y-y'|+|z-z'|).$$

{\rm (iv)}  strong convergence: if $(y_m,z_m)\rightarrow (y,z)$ then
$$\underline{f}_m(t,y_m,z_m)\rightarrow f(t,y,z),
\mbox{and}\ \overline{f}_m(t,y_m,z_m)\rightarrow f(t,y,z),\mbox{as}\
m \rightarrow \infty. $$
\end{lemma}
\begin{lemma}\label{lem:2.4} We assume
$(\underline{Y}^m_{t},\underline{Z}^m_{t})\in S^2\times M^2$ and
$(\overline{Y}^m_{t},\overline{Z}^m_{t})\in S^2\times M^2$ are the
unique solutions of the BDSDEs $(\underline{f}_m,g,T,\xi)$ and
$(\overline{f}_m,g,T,\xi)$ respectively. Then
$$
(\underline{Y}^m_{t},\underline{Z}^m_{t})_{t\in[0,T]}\to
(\underline{Y}_{t},\underline{Z}_{t})_{t\in[0,T]},
$$
and
$$
 (\overline{Y}^m_{t},\overline{Z}^m_{t})_{t\in[0,T]}\to
(\overline{Y}_{t},\overline{Z}_{t})_{t\in[0,T]}, \ (m\to
\infty)
$$
in $S^2 \times M^2,$ where
$(\underline{Y}_{t},\underline{Z}_{t})_{t\in[0,T]}$ and
$(\overline{Y}_{t},\overline{Z}_{t})_{t\in[0,T]}$ are the minimal
solution and maximal solution of BDSDE (\ref{eq:1}).
\end{lemma}

\section{ A simple case: continuous dependence with respect to terminal condition}
\label{sec:3}

This section is devoted to the equivalence of unique solution and
continuous dependence with respect to terminal value $\xi$. Our main result is:
\begin{theorem}\label{thm:3.1}
If Assume (H1)-(H3) hold
for $f$ and $g$, then the following two statements are equivalent.

(i)\ Uniqueness: The equation (\ref{eq:1}) has a unique solution.

(ii)\ Continuous dependence with respect to $\xi$:\  For any
 $\{\xi_n\}_{n=1}^\infty, \xi \in L^2(\Omega,{\cal{F}}_T$, $P)$,
if $\xi_n\rightarrow \xi$ in $L^2(\Omega,{\cal{F}}_T,P)$ as $n \rightarrow \infty,$ then
\begin{equation}\label{eq:2}
\begin{array}{lll}
\lim_{n \rightarrow \infty}E[\sup\limits_{t\in[0,T]}|Y_t^{\xi_n}-Y_t^{\xi}|^2]=0,
\end{array}
\end{equation}
where $(Y_t^{\xi}, Z_t^{\xi})_{t\in[0,T]}$ is any solution of BDSDEs
(\ref{eq:1}) and $(Y_t^{\xi_n}, Z_t^{\xi_n})_{t\in[0,T]}$ are any solutions
of BDSDEs $(f,g,T,\xi^n).$
\end{theorem}
{\bf Proof.}\quad Firstly, we prove that (i)
implies (ii). Given $n$, we note that for any solution
$(Y_t^{\xi_n},$ $Z_t^{\xi_n})_{t\in[0,T]}$ of BDSDEs $(f, g, T,
\xi_n)$, we have
\begin{equation}\label{eq:3}
\begin{array}{lll}
\underline{Y}_{\ t}^{\xi_n}\leq Y_t^{\xi_n}\leq \overline{Y}_{t}^{\
\xi_n}, \mbox{\ P-a.s.\quad $t\in[0,T]$,}
\end{array}
\end{equation}
where $\underline{Y}_{\ t}^{\xi_n}$ and  $\overline{Y}_{t}^{\
\xi_n}$ are the minimal and maximal solutions of BDSDE $(f, g, T$,
$\xi_n)$, respectively.

Now, we consider the following equations:
\begin{equation}\label{eq:4}
\begin{array}{lll}
\underline{Y}_{\ t}^{m,
\xi_n}&=&\xi_n+\int_t^T\underline{f}_m(s,\underline{Y}_{\ s}^{m,
\xi_n}, \underline{Z}_{\ s}^{m, \xi_n})ds\\
 &&+\int_t^Tg(s,\underline{Y}_{\ s}^{m, \xi_n},\underline{Z}_{\ s}^{m, \xi_n})dB_s
 -\int_t^T\underline{Z}_{\ s}^{m, \xi_n}dW_s
\end{array}
\end{equation}
and
\begin{equation}\label{eq:5}
\begin{array}{lll}
\overline{Y}_{t}^{\ m,
\xi_n}&=&\xi_n+\int_t^T\overline{f}_m(s,\overline{Y}_{s}^{\ m, \xi_n},
\overline{Z}_{s}^{\ m, \xi_n})ds\\
 &&+\int_t^Tg(s,\overline{Y}_{s}^{\ m, \xi_n},\overline{Z}_{s}^{\ m, \xi_n})dB_s
 -\int_t^T\overline{Z}_{s}^{\ m, \xi_n}dW_s
\end{array}
\end{equation}
where $(\underline{Y}_{\ t}^{m, \xi_n},\underline{Z}_{\ t}^{m, \xi_n})_{t\in[0,T]}$ and
 $(\overline{Y}_{t}^{\ m, \xi_n},\overline{Z}_{t}^{\ m, \xi_n})_{t\in[0,T]}$ are unique solutions of (\ref{eq:4})
and (\ref{eq:5}) respectively.

Thanks to Lemma \ref{lem:2.4}, we know that
$$(\underline{Y}_{\ t}^{m, \xi_n},\underline{Z}_{\ t}^{m, \xi_n})
\rightarrow (\underline{Y}_{\ t}^{\xi_n},\underline{Z}_{\
t}^{\xi_n})\ \mbox{and}\ (\overline{Y}_{t}^{\ m,
\xi_n},\overline{Z}_{t}^{\ m, \xi_n}) \rightarrow
(\overline{Y}_{t}^{\ \xi_n},\overline{Z}_{t}^{\ \xi_n}),\
t\in[0,T]$$ in $S^2\times M^2$ as $m\rightarrow \infty,$ and from
Comparison Theorem 3.1 of \cite{SZ} get the following inequalities
\begin{equation}\label{eq:6}
\begin{array}{lll}
\underline{Y}_{\ t}^{m, \xi_n}\leq\underline{Y}_{\ t}^{\xi_n}\leq
Y_t^{\xi_n} \leq \overline{Y}_{t}^{\ \xi_n}\leq\overline{Y}_{t}^{\
m, \xi_n}, \mbox{\ for\ any\ $n, t\in[0,T]$\ and\ m$\geq$ K}.
\end{array}
\end{equation}
From (\ref{eq:6}), we have
\begin{eqnarray*}
Y_{t}^{\xi_n}-Y_{t}^{\xi} &=&Y_{t}^{\xi_n}-\overline{Y}_{t}^{\
m,\xi_n}+\overline{Y}_{t}^{\ m,\xi_n}
-\overline{Y}_{t}^{\ m,\xi}+\overline{Y}_{t}^{\ m,\xi}-Y_{t}^{\xi}\\
&\leq& (\overline{Y}_{t}^{\ m,\xi_n} -\overline{Y}_{t}^{\
m,\xi})+(\overline{Y}_{t}^{\ m,\xi}-Y_{t}^{\xi}),
\end{eqnarray*}
and
\begin{eqnarray*}
Y_{t}^{\xi_n}-Y_{t}^{\xi} &=&Y_{t}^{\xi_n}-\underline{Y}_{\
t}^{m,\xi_n}+\underline{Y}_{\ t}^{m,\xi_n}
-\underline{Y}_{\ t}^{m,\xi}+\underline{Y}_{\ t}^{m,\xi}-Y_{t}^{\xi}\\
&\geq& (\underline{Y}_{\ t}^{m,\xi_n} -\underline{Y}_{\
t}^{m,\xi})+(\underline{Y}_{\ t}^{m,\xi}-Y_{t}^{\xi}).
\end{eqnarray*}
Thus
\begin{eqnarray*}
&&E[\sup\limits_{t\in[0,T]}|Y_t^{\xi_n}-Y_t^{\xi}|^2]\\
&\leq&2E[\sup\limits_{t\in[0,T]}|\overline{Y}_{t}^{\
m,\xi_n}-\overline{Y}_{t}^{\ m,\xi}|^2]
+2E[\sup\limits_{t\in[0,T]}|\overline{Y}_{t}^{\ m,\xi}-Y_{t}^{\xi}|^2]\\
&+&2E[\sup\limits_{t\in[0,T]}|\underline{Y}_{t}^{\
m,\xi_n}-\underline{Y}_{t}^{\ m,\xi}|^2] +2E[\sup\limits_{t\in[0,T]}|\underline{Y}_{t}^{\
m,\xi}-Y_{t}^{\xi}|^2],
\end{eqnarray*}
where $(\underline{Y}_{\ t}^{m, \xi_n},\underline{Z}_{\ t}^{m, \xi_n})_{t\in[0,T]}$ and
 $(\overline{Y}_{t}^{\ m, \xi_n},\overline{Z}_{t}^{\ m, \xi_n})_{t\in[0,T]}$ are unique
solutions of\\
BDSDEs $(\underline{f}_m,$ $g,T,\xi)$ and
$(\overline{f}_m,g,T,\xi)$ respectively.

By Lemma \ref{lem:2.2} and Lemma \ref{lem:2.3}, as $n\rightarrow \infty$, we have
$$E[\sup\limits_{t\in[0,T]}|\underline{Y}_{\ t}^{m,\xi_n}
-\underline{Y}_{\ t}^{m,\xi}|^2]\rightarrow 0,
\mbox{\ and\ }E[\sup\limits_{t\in[0,T]}|\overline{Y}_{\ t}^{\ m,\xi_n}
-\overline{Y}_{\ t}^{\ m,\xi}|^2]\rightarrow 0, \mbox{\ for\ any\
}m.$$

By Lemma \ref{lem:2.4} and the uniqueness of solution for BDSDEs (\ref{eq:1}), we get
$$E[\sup\limits_{t\in[0,T]}|\underline{Y}_{\ t}^{m,\xi}
-\underline{Y}_{\ t}^{\xi}|^2]\rightarrow 0 \mbox{\ and\ }
E[\sup\limits_{t\in[0,T]}|\overline{Y}_{t}^{\ m,\xi}
-\overline{Y}_{t}^{\ \xi}|^2]\rightarrow 0$$ as $m\rightarrow
\infty.$ That is (ii).

Now, we prove that (ii) implies (i). We take $\xi_n=\xi.$ For
equation $(f,g,T,\xi_n),$ we set $Y_t^{\xi_n}=\overline{Y}_{t}^{\
\xi_n}=\overline{Y}_{t}^{\ \xi}.$ For the equation (\ref{eq:1}), we set
$Y_t^{\xi}=\underline{Y}_{\ t}^{\xi}.$ For (ii), we have
$\overline{Y}_{t}^{\ \xi}=\underline{Y}_{\ t}^{\xi}.$ \quad$\Box$

\begin{remark}\label{rmk:3.2}
In fact, when the solution of (\ref{eq:1}) is not unique, the continuous
dependence may not hold true in general. For example, we take
$f(t,y,z)=3y^{2/3}$, $\xi= 0$ and $g$ such that $g(t,y,0)=0$
for all $t\in [0,T]$, $(y,z)\in \mathbb{R}\times  \mathbb{R}^d$.
It is easy to know that $(y_t,z_t)_{t\in [0,T]} = (0, 0)_{t\in [0,T]}$
and $(Y_t,Z_t)_{t\in [0,T]} = ((T-t)^3, 0)_{t\in [0,T]}$ both are solutions of BDSDE
\begin{equation*}
Y_t=\int_t^T3Y_s^{2/3}ds +\int_t^Tg(s,Y_s,Z_s)dB_s-\int_t^TZ_sdW_s,
\quad 0\leq t\leq T.
\end{equation*}
Set $\xi_n = 1/n$, the BDSDEs
\begin{equation*}
Y_t=\displaystyle\frac{1}{n}+ \int_t^T3Y_s^{2/3}ds
+\int_t^Tg(s,Y_s,Z_s)dB_s-\int_t^TZ_sdW_s,
\ 0\leq t\leq T,\ n = 1, 2,\cdots
\end{equation*}
have unique solutions $(y^{\frac{1}{n}}_t, z^{\frac{1}{n}}_t)
=((T-t+\frac{1}{\sqrt{n}})^3, 0)$ for $n = 1, 2,\cdots $. But
$$\lim\limits_{n\to \infty}E[ \sup\limits_{t\in [0,T]}
|y^{\frac{1}{n}}_t-y_t|^2] = T^6 \not= 0 =
\lim\limits_{n\to \infty}E[ \sup\limits_{t\in [0,T]}
|y^{\frac{1}{n}}_t-Y_t|^2].$$
\end{remark}

\section{The general case}\label{sec:4}

In this section, we will deal with the more general case, that is, the relationship
between uniqueness of solution and the continuous dependence with respect not
only to $\xi$ but also to $f$ and $g$. Now, we consider the following BDSDEs:
\begin{equation}\label{eq:7}
\begin{array}{lll}
Y^{\lambda}_t=\xi^{\lambda}+\int_t^Tf^{\lambda}(s,Y^{\lambda}_s,Z^{\lambda}_s)ds +\int_t^Tg^{\lambda}(s,Y^{\lambda}_s,Z^{\lambda}_s)dB_s-\int_t^TZ^{\lambda}_sdW_s,
\end{array}
\end{equation}
where $\lambda$ belongs to a nonempty set $D\subset R^n.$ The coefficients
$$f^{\lambda}(\omega,t,y,z):\Omega\times [0,T]\times \mathbb{R}
\times \mathbb{R}^{ d}\rightarrow \mathbb{R},
\mbox{\ and\ } g^{\lambda}(\omega,t,y,z):\Omega\times [0,T]\times
\mathbb{R}\times \mathbb{R}^{ d}\rightarrow \mathbb{R}^l,$$ satisfying the following
conditions:

\vspace{0.1cm}(H1') linear growth: $\exists\ 0<K< \infty$, such that

\vspace{0.1cm}$|f^\lambda(\omega,t,y,z)|\leq K(1+|y|+|z|), \quad\forall\ \lambda,
\omega, t, y, z\in D\times\Omega \times [0,T]\times \mathbb{R} \times \mathbb{R}^d.$

\vspace{0.1cm}(H2') For fixed $\lambda,\ \omega$ and $t$, $f^\lambda(\omega, t, \cdot,\cdot)$ is
continuous.

(H3')  uniform continuity:  $f^{\lambda}$ and  $g^{\lambda}$ are
continuous in $\lambda=\lambda_0$ uniformly with respect to $(y,z).$

(H4') For fixed $\lambda$, there exist constants $c>0$ and
$0<\alpha<1$ such that

\vspace{0.1cm}\hspace{2cm}$|g^{\lambda}(\omega,t,y^1,z^1)
-g^{\lambda}(\omega,t,y^2,z^2)|^2\leq c|y^1-y^2|^2+\alpha
|z^1-z^2|^2,$

\vspace{0.1cm}\noindent for all $(\omega,t)\in \Omega\times [0,T],
\ (y^1,z^1)\in R\times \mathbb{R}^d, \ (y^2,z^2)\in \mathbb{R}\times \mathbb{R}^d.$

(H5') Lipschitz condition: $\exists\ 0<c< \infty$, such that

\vspace{0.1cm}\hspace{2cm}$|f^{\lambda}(\omega,t,y^1,z^1)
-f^{\lambda}(\omega,t,y^2,z^2)|^2\leq c(|y^1-y^2|^2+|z^1-z^2|^2),$

\vspace{0.1cm}\noindent for all $(\omega,t)\in \Omega\times [0,T],
\ (y^1,z^1)\in \mathbb{R}\times \mathbb{R}^d, \ (y^2,z^2)\in \mathbb{R}\times \mathbb{R}^d.$

Under (H3')-(H5'), the BDSDE (\ref{eq:7}) has a unique adapted solution for any $\lambda\in D.$
And we have the following property:

\begin{lemma}\label{lem:4.1}
If $\xi^{\lambda}\rightarrow \xi^{\lambda_0}$ in
$L^2(\Omega,{\cal{F}}_T,P)$ as $\lambda \rightarrow \lambda_0,$
assumptions (H3')-(H5') hold for $f^{\lambda}$ and $g^{\lambda}$.
Moreover $(Y^{\lambda}_t, Z^{\lambda}_t)_{t\in[0,T]}$ and
$(Y^{\lambda_0}_t; Z^{\lambda_0}_t)_{t\in[0,T]}$ are the solutions
of the BDSDEs $(f^{\lambda},g^{\lambda},T,\xi^{\lambda})$ and
$(f^{\lambda_0},g^{\lambda_0},T,\xi^{\lambda_0})$ respectively, then
\begin{equation}\label{eq:8}
\begin{array}{lll}
&&E[\sup\limits_{t\in[0,T]}|Y^{\lambda}_t-Y^{\lambda_0}_t|^2]\\
&\leq&
CE|\xi^{\lambda}-\xi^{\lambda_0}|^2 +CE\int_0^T|f^{\lambda}(t,Y^{\lambda}_t, z^{\lambda}_t)
-f^{\lambda_0}(t,Y^{\lambda_0}_t, z^{\lambda_0}_t)|^2dt\\&&+CE\int_0^T|g^{\lambda}(t,Y^{\lambda}_t, z^{\lambda}_t)
-g^{\lambda_0}(t,Y^{\lambda_0}_t, z^{\lambda_0}_t)|^2dt,
\end{array}
\end{equation}
where $C$ is a positive constant only depending on Lipschitz constan
t $c$ and $\alpha$. Moreover, we have
\begin{equation}\label{eq:9}
\begin{array}{lll}
\lim\limits_{\lambda \rightarrow \lambda_0}E[\sup\limits_{t\in[0,T]}
|Y^{\lambda}_t-Y^{\lambda_0}_t|^2]=0.
\end{array}
\end{equation}
\end{lemma}
{\bf Proof.}\quad By the usual techniques
of BDSDEs we can get inequality (\ref{eq:8}) (see \cite{PP2} for detail).
Because of the continuity of $f^\lambda$ and  $g^\lambda$ in $\lambda =
\lambda_0$ and Lebesgue dominated convergence theorem we take limit
both sides of (\ref{eq:8}) and get equation (\ref{eq:9}). \quad$\Box$

Now, we introduce the approximation sequences of $f^\lambda$ as follows:
\begin{eqnarray*}
\underline{f}^{\ \lambda}_{\ m}(\omega,t,y,z)=\inf_{(y', z')\in
Q^{1+d}}
\{f^\lambda(\omega,t,y',z')+m(|y-y'|+|z-z'|)\} 
\end{eqnarray*}
and
\begin{eqnarray*}
\overline{f}^{\ \lambda}_m(\omega,t,y,z)=\sup_{(y', z')\in Q^{1+d}}
\{f^\lambda(\omega,t,y',z')-m(|y-y'|+|z-z'|)\} 
\end{eqnarray*}

\begin{lemma}\label{lem:4.2}  If $f^\lambda$ satisfies
Assumptions (H1')-(H3'), then for any $m\geq K$, we have

(i)  linear growth: $\forall\ ( y, z)\in \mathbb{R} \times \mathbb{R}^d$ and $t\in
[0,T],$ $$|\underline{f}^\lambda_{\ m}(t,y,z)|\leq K(1+|y|+|z|),
\mbox{\ and\ } |\overline{f}^{\ \lambda}_m(t,y,z)|\leq
K(1+|y|+|z|).$$

(ii)  monotonicity in $m$: $\forall\ ( y, z)\in \mathbb{R} \times \mathbb{R}^d$ and
$t\in [0,T],$ $\underline{f}^\lambda_{\ m}(t,y,z)$ is non-decreasing
in $m$ and $\overline{f}^{\ \lambda}_m(t,y,z)$ is non-increasing in
$m$.

(iii) Lipschitz condition: $\forall\ y_1, y_2 \in \mathbb{R}, z_1,z_1\in \mathbb{R}^d$
and $t\in [0,T],$
$$|\underline{f}^\lambda_{\ m}(t,y,z)-\underline{f}^\lambda_m(t,y',z')|\leq
m(|y-y'|+|z-z'|),$$ and $$|\overline{f}^{\
\lambda}_m(t,y,z)-\overline{f}^{\ \lambda}_m(t,y',z')|\leq
m(|y-y'|+|z-z'|).$$

(iv)  strong convergence: if $(y_m,z_m)\rightarrow (y,z)$ as $m
\rightarrow \infty,$ then
$$\underline{f}^\lambda_{\ m}(t,y_m,z_m)\rightarrow f^\lambda(t,y,z),
\mbox{\ and\ } \overline{f}^{\ \lambda}_m(t,y_m,z_m)\rightarrow
f^\lambda(t,y,z)\mbox{\ as\ } m \rightarrow \infty. $$

(v) Both $\underline{f}^\lambda_{\ m}$ and $\overline{f}^{\
\lambda}_m$ are continuous in $\lambda=\lambda_0.$
\end{lemma}

\noindent{\bf Proof.}\quad It is easy to check (i)-(iv)
(see \cite{SGL}). Now, we prove (v). For any $\varepsilon > 0$, by the
definition of $\underline{f}^\lambda_{\ m},$ there exist
$(y^{\varepsilon,\lambda}, z^{\varepsilon,\lambda})$ and
$(y^{\varepsilon,\lambda_0}, z^{\varepsilon,\lambda_0})$ such that
\begin{eqnarray*}
&&f^\lambda(t,y^{\varepsilon,\lambda}, z^{\varepsilon,\lambda})
+m(|y-y^{\varepsilon,\lambda}|+|z-z^{\varepsilon,\lambda}|)-\varepsilon
\leq \underline{f}^\lambda_{\ m}(t,y,z)\\
&&\leq f^\lambda(t,y^{\varepsilon,\lambda_0},
z^{\varepsilon,\lambda_0})
+m(|y-y^{\varepsilon,\lambda_0}|+|z-z^{\varepsilon,\lambda_0}|),
\end{eqnarray*}
and
\begin{eqnarray*}
&&f^{\lambda_0}(t,y^{\varepsilon,\lambda_0}, z^{\varepsilon,\lambda_0})
+m(|y-y^{\varepsilon,\lambda_0}|+|z-z^{\varepsilon,\lambda_0}|)-\varepsilon
\leq \underline{f}^{\lambda_0}_{\ m}(t,y,z)\\
&&\leq f^{\lambda_0}(t,y^{\varepsilon,\lambda},
z^{\varepsilon,\lambda})
+m(|y-y^{\varepsilon,\lambda}|+|z-z^{\varepsilon,\lambda}|),
\end{eqnarray*}
thus
\begin{eqnarray*}
&&f^\lambda(t,y^{\varepsilon,\lambda}, z^{\varepsilon,\lambda})
-f^{\lambda_0}(t,y^{\varepsilon,\lambda}, z^{\varepsilon,\lambda})-\varepsilon
\leq \underline{f}^\lambda_{\ m}(t,y,z)-\underline{f}^{\lambda_0}_{\ m}(t,y,z)\\
&&\leq f^\lambda(t,y^{\varepsilon,\lambda_0},
z^{\varepsilon,\lambda_0})
-f^{\lambda_0}(t,y^{\varepsilon,\lambda_0},
z^{\varepsilon,\lambda_0})+\varepsilon.
\end{eqnarray*}
Because $f^\lambda$ is continuous in $\lambda=\lambda_0$ uniformly
with respect to $(y,z)$, we obtain the continuity of
$\underline{f}^\lambda_{\ m}$ and $\overline{f}^{\ \lambda}_m$ in
$\lambda=\lambda_0$. \quad$\Box$

\begin{lemma}\label{lem:4.3}
If $f^\lambda$ and $g^\lambda$ satisfy (H1')-(H4'),
 and the processes $(\underline{Y}^{\lambda,m}_t,\underline{Z}^{\lambda,m}_t)_{t\in[0,T]}$ and
 $(\overline{Y}^{\lambda,m}_t,$ $\overline{Z}^{\lambda,m}_t)_{t\in[0,T]}$
are the unique solutions of the BDSDEs $(\underline{f}^{\lambda,m},g^\lambda,T,\xi^\lambda)$ and $(\overline{f}^{\lambda,m},g^\lambda,T,\xi^\lambda)$ respectively, then, for any $\lambda\in
D$, we have
$$
(\underline{Y}^{\lambda,m}_t,
\underline{Z}^{\lambda,m}_t)_{t\in[0,T]}
\rightarrow(\underline{Y}^{\lambda}_t,
\underline{Z}^{\lambda}_t)_{t\in[0,T]},$$
and
$$(\overline{Y}^{\lambda,m}_t, \overline{Z}^{\lambda,m}_t)_{t\in[0,T]} \rightarrow(\overline{Y}^{\lambda}_t,
\overline{Z}^{\lambda}_t)_{t\in[0,T]}$$
in $S^2\times M^2$ as $m\rightarrow\infty,$ where
$(\underline{Y}^{\lambda}_t, \underline{Z}^{\lambda}_t)_{t\in[0,T]}$
and
$(\overline{Y}^{\lambda}_t, \overline{Z}^{\lambda}_t)_{t\in[0,T]}$
are the minimal solution and maximal solution of BDSDE {\rm (\ref{eq:7})}.
\end{lemma}

\noindent Now, we give our result for the general case.

\begin{theorem}\label{thm:4.4}
If $f^\lambda$ and $g^\lambda$ satisfy {\rm (H1')-(H4')}, then the following
statements are equivalent.

{\rm (iii)} Uniqueness: there exists a unique solution of BDSDE {\rm (\ref{eq:7})} with
$\lambda=\lambda_0,$ that is, the solution of
$(f^{\lambda_0},g^{\lambda_0},T,\xi^{\lambda_0})$  is unique.

{\rm (iv)} Continuous dependence with respect to $f$, $g$ and $\xi$: for any
$\xi^{\lambda},\xi^{\lambda_0} \in L^2(\Omega,{\cal{F}}_T,P)$, if
$\xi^{\lambda}\rightarrow\xi^{\lambda_0}$ in
$L^2(\Omega,{\cal{F}}_T,P)$ as $\lambda\rightarrow \lambda_0$,
$(Y^{\lambda}_t,Z^{\lambda}_t)_{t\in[0,T]}$are any solutions of
BDSDE {\rm (\ref{eq:7})},  $(Y^{\lambda_0}_t,Z^{\lambda_0}_t)_{t\in[0,T]}$ is any
solution of BDSDE {\rm (\ref{eq:7})} with $\lambda = \lambda_0$, then
$$\lim\limits_{\lambda \rightarrow \lambda_0}E[\sup\limits_{t\in[0,T]}
|Y^{\lambda}_t-Y^{\lambda_0}_t|^2]=0.$$
\end{theorem}

\noindent{\bf Proof.}\ This proof is similar to that
of Theorem 3.1. For the sake of completeness, we give the sketch of
proof. Firstly, we prove (iii) implies (iv). We can get the
inequalities similarly to (\ref{eq:6}), that is, from Comparison Theorem 3.1
of \cite{SZ} get the following inequalities
\begin{eqnarray*}
\underline{Y}_{\ t}^{\lambda,m}\leq\underline{Y}_{\ t}^{\
\lambda}\leq Y_t^{\lambda} \leq \overline{Y}_{t}^{\
\lambda}\leq\overline{Y}_{t}^{\lambda,m}, \mbox{\ for\ any
$t\in[0,T]$\ and\ m$\geq$ K.}
\end{eqnarray*}
So
\begin{eqnarray*}
&&E[\sup\limits_{t\in[0,T]}|Y_t^{\lambda}-Y_t^{\lambda_0}|^2]\\
&\leq&2E[\sup\limits_{t\in[0,T]}|\underline{Y}_{\
t}^{\lambda,m}-\underline{Y}_{\ t}^{\lambda_0,m}|^2]
+2E[\sup\limits_{t\in[0,T]}|\underline{Y}_{\ t}^{\lambda_0,m}-Y_{t}^{\lambda_0}|^2]\\
&+&2E[\sup\limits_{t\in[0,T]}|\overline{Y}_{t}^{\lambda,m}-\overline{Y}_{t}^{\lambda_0,m}|^2] +2E[\sup\limits_{t\in[0,T]}|\overline{Y}_{t}^{\lambda_0,m}-Y_{t}^{\lambda_0}|^2].
\end{eqnarray*}
Fixed $m$, by Lemma \ref{lem:4.1} and Lemma \ref{lem:4.2}, and the continuity of
$\underline{f}_{\ m}^\lambda$ and $\overline{f}_m^\lambda$ in
$\lambda=\lambda_0,$ we have
$$E[\sup\limits_{t\in[0,T]}|\underline{Y}_{\  t}^{\lambda,m}
-\underline{Y}_{\  t}^{\lambda_0,m}|^2]\rightarrow 0 \mbox{\ and\
} E[\sup\limits_{t\in[0,T]}|\overline{Y}_{t}^{\  \lambda,m}
-\overline{Y}_{t}^{\  \lambda_0,m}|^2]\rightarrow 0$$ as
$\lambda\rightarrow\lambda_0,$ for any $m\geq K.$ By Lemma \ref{lem:4.3} and the
uniqueness of solution for BDSDEs
$(f^{\lambda_0},g^{\lambda_0},T,\xi^{\lambda_0})$ (Condition (iii)),
we get , as $m\rightarrow \infty,$
$$E[\sup\limits_{t\in[0,T]}|\underline{Y}_{\  t}^{\lambda_0,m}
-\underline{Y}_{\  t}^{\lambda_0}|^2]\rightarrow 0 \mbox{\ and\ }
E[\sup\limits_{t\in[0,T]}|\overline{Y}_{t}^{\  \lambda_0,m}
-\overline{Y}_{t}^{\ \lambda_0}|^2]\rightarrow 0.$$
 \noindent This implies (iv).

Now, we prove that (iv) implies (iii). We take
$\xi^{\lambda}=\xi^{\lambda_0},f^{\lambda}=f^{\lambda_0},g^{\lambda}=g^{\lambda_0}.$
For equation (\ref{eq:7}), set
$Y_t^{\lambda}=\overline{Y}_{t}^{\lambda}=\overline{Y}_{t}^{\lambda_0}.$
For the equation $(f^{\lambda_0},g^{\lambda_0},T,\xi^{\lambda_0})$,
we set $Y_t^{\lambda_0}=\underline{Y}_{\ t}^{\lambda_0}.$ For (iv),
we have $\overline{Y}_{t}^{\lambda_0}=\underline{Y}_{\
t}^{\lambda_0}.$ \quad$\Box$

\begin{remark}\label{rmk:4.5}
In the standard situation where $f$ satisfies linear growth condition and Lipschitz
condition in $(y,z)$, it has been proved by Pardoux and Peng \cite{PP2} that there exists
a unique solution. In this case, the continuous dependence with respect to $f$, $g$
and $\xi$ is is described by the  inequality (\ref{eq:8}) (see  \cite{SGL}).
Our result in this paper, which can be regarded as the analog of the inequality
(\ref{eq:8}) in some sense, provides a useful method to study BSDEs with continuous coefficient.
\end{remark}

\end{document}